\documentstyle[11pt] {article}

\begin{document}
\def \om {\Omega}
\def \ra {\rightarrow}
\def \o  {\overline}
\def \p  {\partial}
\def \l  {\lambda}
\def\endpf{\hbox{\vrule height1.5ex width.5em}}
\def \d  {\displaystyle}
\def \f  {\frac}
\def \e  {\varepsilon}
\def \un {\underline{\lim}}

\date{}
\title{Decomposition of vector-valued  divergence free Sobolev functions  and shape optimization for stationary Navier-Stokes
equations
   \footnote{
This work was supported by  National Natural Science Foundation of
China, Grant No.  10471053, 60574017 and by the key programme of
the Ministry of National Education of China.}}
\author{Gengsheng Wang \\
 School of Mathematics, Wuhan University,\\ Wuhan, Hubei,
430072, P.R.of China\\
\verb "wanggs@public.wh.hb.cn"\\
\ Donghui Yang\\
 School of Mathematics, Wuhan University,\\ Wuhan, Hubei,
430072, P.R.of China
\\ \verb "dongfyang@yahoo.com.cn"} \maketitle

\begin{abstract}\quad We establish a divergence free partition for
vector-valued Sobolev functions with free divergence in
 ${\bf R}^n, n\geq 1$. We  prove that  for any domain $\om$ of
 class $\cal C$ in ${\bf R}^n,n=2,3$, the space
$D_0^1(\om)\equiv\{{\mathbf{v}} \in H^1_0(\Omega)^n ;
\mbox{div}{\mathbf{v}}=0\}$ and the space $H_{0,\sigma}^1(\om)\equiv
\overline{\{{\mathbf{v}}\in
C^{\infty}_0(\Omega)^n;\mbox{div}{\mathbf
v}=0\}}^{\|\cdot\|_{H^1(\om)^n}}$, which is  the completion of
$\{{\mathbf{v}} \in C^{\infty}_0(\Omega)^n; \mbox{div}{\mathbf
v}=0\}$ in the  $H^1(\Omega)^n$-norm, are identical. We will also prove that
$H_{0,\sigma}^1(D\setminus\overline\om)=\{{\mathbf v}\in
H_{0,\sigma}^1(D); {\mathbf v}=0\ \mbox{ a.e. in } \om\}$,  where
$D$ is a  bounded  Lipschitz domain such that
$\om\subset\subset D$. These results, together with  properties for
domains of class $\mathcal C$, are used to  solve an  existence problem  in the shape
optimization theory  of the stationary Navier-Stokes equations.

 \vskip 0.1cm


\vskip 0.1cm


\end{abstract}

\def\theequation{1.\arabic{equation}}
\setcounter{equation}{0} \noindent {\large\bf1.\quad Introduction
}

\vskip 3mm

  In this paper, we study an existence problem in the shape optimization
  theory  of the stationary Navier-Stokes equation over
  $D \backslash \overline\Omega$ with $\Omega\subset\subset D$ varying. Here,
  $D\subset {\bf R}^n$, $n=2, 3$, is a given
   bounded Lipschitz domain and $\Omega$ is  supposed to satisfy the property
   $\mathcal C$ to  be defined  in section 2.
  Our main   ingredient for such an investigation is  a new divergence-free partition for
  divergence free vector-valued Sobolev functions defined over
  $D \backslash \overline\Omega$. This partition is established by making use of
    the  classical  Hodge Theory in Differential Geometry.
  To be more precise, we let  $D$ be an open subset  in ${\bf R}^n$
  ($n\geq 1$), and let $\{U_j\}_{j=1}^m $ be an open covering of $\overline D$. Suppose that
  ${\mathbf
  u} \in D^1_0(D) \equiv\{{\mathbf u }\in H_0^1(D)^n; div{\mathbf u}=0\}
  $. Then  our first main result  states that
  there are functions $\{{\mathbf u }_j\}_{j=1}^m$ such that
  ${\mathbf u}_j\in D_0^1(U_j)$, ${\mathbf u }=\sum_{j=1}^m {\mathbf
  u_j}$ over $D$ and $\|{\mathbf u}_j\|_{H^1(U_j)^n} \leq
  C \|{ \mathbf u }\|_{H^1(D)^n}$ for some positive constant C
  independent of ${\mathbf{u}}$.  Moreover,
  such a partition inherits the following local property:
 If ${\mathbf u}\!\!\mid_\om =0$ with
  $\om \subset \subset D$  and
$\om $ is of class $\mathcal C$ to be defined in section 2, then
we can choose $\{{\mathbf u}_j\}_{j=1}^m$ such that ${\mathbf
u}_j\!\!\mid_\om=0$ for  $j=1,\cdots,m$.

The above decomposition, besides its application to the shape
optimization problem for the stationary Navier-Stokes equation,
has many corollaries which might be interesting in their own
right. For instance, it gives the following property for the
function spaces $H_{0,\sigma}^1(\om)$ and $H_{0,\sigma}^1(D
\backslash\overline{\om})$ where
$$H_{0,\sigma}^1(\om)\equiv\overline{\{ {\mathbf u}\in
C_0^\infty(\om)^n; \mbox{div}{\mathbf u}=0}\}^{\| \cdot
\|_{H^1(\om)^n}}$$ (the completion of $ \{ {\mathbf u}\in
C_0^\infty(\om)^n; \mbox{div}{\mathbf u}=0 \}$ in the norm of
$H^1(\om)^n$):

If $\om$ is of class $\cal C$ and $D$ is Lipschitz,
then
\begin{eqnarray}\label{1.1}
 && H^1_{0,\sigma}(\om)=D^1_0(\om),\ \\ \label{1.2}  && H_{0,\sigma}^1(D\setminus\overline\om)=\{{\mathbf v}\in
H^1_{0,\sigma}(D); {\mathbf v}=0 \ \ \mbox{a.e.\ in }\ \om \}.
\end{eqnarray}

\noindent
(1.2)  has   the following consequences:

(i)\quad If ${\mathbf v}\in
H_{0,\sigma}^1(D\setminus\overline\om)$, then $\tilde {\mathbf
v}\in H_{0,\sigma}^1(D)$, where
\begin{eqnarray*}
\tilde {\mathbf v} =\begin{cases}{
{\mathbf v}\ \ &{\rm in}\ \
$D\setminus\overline\om$,\cr 0,\ \ &{\rm in }\ \ $\om$.\cr}
\end{cases}
\end{eqnarray*}

(ii)\quad If ${\mathbf v}\in H_{0,\sigma}^1(D)$ and $ {\mathbf
v}=0$ a.e. in $\om$, then $\hat {\mathbf v}\in
H_{0,\sigma}^1(D\setminus\overline\om)$, where $\hat {\mathbf v}$
is the restriction of ${\mathbf v}$ over $D\setminus\overline\om$.

 (1.1) and (1.2) can  be compared with the classical
result in $\cite{kn:[12]}$: If $v\in H_0^1(D) $ and $v=0$
quasi-everywhere in $D\setminus\om$ ( where $\om$ is an open set
), then $v\in H_0^1(\om)$. They can also be compared with  the
result in $\cite{kn:[9]}$, which states that if $v\in H_0^1(D)$
and $v=0$ a.e. in $D\setminus\om$ ( where $\om$ is of class $\cal
C$ ), then $v\in H_0^1(\om)$. We notice that  (1.1)-(1.2) have
been obtained earlier in $\cite{kn:[4]}$ and $\cite{kn:[7]}$ under
the assumption that $\om$ is Lipschitz.
 The novelty in (1.1)-(1.2) is that $\om$ does not have to
 have a Lipschitz boundary. In fact,  a generic domain   with property
$\cal C$ is not a finite union
of  Lipschitz domains. It may have cusps and infinitely many
oscillations with vanishing amplitude.

Fundamentally based on the aforementioned  decomposition  and
properties on the family $\mathcal O$ which will be defined in
Section 2, we study in section 5 a shape optimization problem. In
this problem, we will look for a domain of class ${\mathcal C}$
which minimizes a certain cost functional associated with the
solutions of
 the stationary Navier-Stokes equation.
 Our main purpose for this part of the paper is to prove Theorem 5.1 to be stated in section 5.

There have been extensive studies  for the  existence problems in the shape
optimization theory for elliptic equations in the past twenty years.  Here,
 we mention the
work in  $\cite{kn:[3]}$, $\cite{kn:[6]}$,
 where results are obtained over domains with sufficiently smooth
boundary. We mention  the work in  $\cite{kn:[8]}$ and $\cite{kn:[2]}$,
which is based on  a certain capacity constraint and penalty terms, respectively.
For more recent work, we mention $\cite{kn:[9]}$ and $\cite{kn:[10]}$
where the domains are assumed to have  a certain segment property.
These two papers are based on the
analysis of the set of bounded open sets of class $\mathcal C$ in the sense
of Maz'ya $\cite{kn:[5]}$, which, in terms of  Adams $\cite{kn:[1]}$, is
identical to the set of open subsets with a certain segment
property.

Our paper is organized as follows. In Section 2, we give some
definitions and preliminary facts  for domains of class $\cal C$. In
Section 3, we establish a divergence free partition of
vector-valued Sobolev functions  mentioned above. In Section 4, we study
properties of the function spaces $H_{0,\sigma}^1(\om)$ and
$H_{0,\sigma}^1(D\setminus\overline\om)$. Finally,  in Section 5, we
prove the existence theorem of a shape optimization problem  for the stationary
Navier-Stokes equation.

\vskip 3mm

\def\theequation{2.\arabic{equation}}
\setcounter{equation}{0} \noindent {\large\bf2.\quad Preliminary facts
for open subsets of class  ${\cal C}$}

\vskip 3mm

  Let $\om$ be a bounded  open subset in ${\bf R}^n,\ n\geq1$.  We say
that $\om $ is of class $\cal C$ or has property ${\cal C}$,  if
the following properties (2.1), (2.2) and (2.3) are satisfied:

 There is a family ${\cal
F}_\Omega $ of real valued continuous functions $g$ defined over
$\o{S(0, k_\om)}$, where  $S(0,k_\om)\subset
{\bf R}^{n-1}$ is  the
open ball  centered at 0 and of radius $k_\om>0$, such that
\begin{eqnarray}\label{2.1}
\p\om=\bigcup_{g\in {\cal F}_\om}\{o_g+R_g((s,g(s))); s\in
S(0,k_\om)\}.\end{eqnarray} Here $o_g\in \p\om$ is a certain point
which gives the center of the local system of axes
 and $R_g: {\bf R}^n\rightarrow {\bf R}^n$ is a certain
unitary rotation  such that $R_g(0,\cdots,0,1)=y_g$, where the
unit vector $y_g$ gives the local vertical direction. Moreover,
for any
 $g\in {\cal F}_\om$, there is an $a_\om>0$ such that
\begin{eqnarray}\label{2.2}
&& o_g+R_g((s,g(s)))-ty_g\in \om, t\in(0,a_\om),s\in
S(0,k_\om),\\
\label{2.2}&& o_g+R_g((s,g(s)))+ty_g\in {\bf R}^n\backslash
\overline\om, t\in(0, a_\om), s\in S(0,k_\om).\end{eqnarray}

 By the Lebesgue Lemma, it  follows from (2.1) that there is an $r_\om\in
(0,k_\om)$ such that the  restricted local charts defined on
$\overline{S(0,r_\om)}$ also give a covering of $\p \om $:
\begin{eqnarray}
\p\om=\bigcup_{g\in{\cal F}_\om}\{o_g+R_g((s,0))+g(s)y_g; s\in
\overline {S(0,r_\om)}\}.\end{eqnarray}

 The above definition was given
in $\cite{kn:[9]}$, which is a slight modification of the
corresponding definition given by Maz'ya $\cite{kn:[5]}$. Roughly
speaking, $\om$ is of class $\cal C$ if $\p\om$ is locally a graph
of a continuous function with segment property. We easily see that
if $\om$ is of class $\cal C$, then there are open sets
$\{V_j\}_{j=1}^m$ such that
\begin{eqnarray} && \p \om\subset \bigcup_{j=1}^m V_j, \\ &&
V_j=\{x\in {\bf R}^n; x=p_j(s)+t{\mathbf v}_j, s\in S(0,k_\om),
t\in (-d, d) \},\end{eqnarray}
 with ${\mathbf v}_j$ a unit vector, $p_j(s)=o_j+R_j((s,g_j(s))),
 g_j\in {\cal F}_\Omega, o_j\in\p\Omega$ and $R_j$ a unitary rotation
 satisfying $R_j((0,\cdots,0,1))={\mathbf v}_j$.
Moreover, by the
segment property in (2.2) and (2.3), we may take $d>0$ small enough
such  that for any $j=1,\cdots,m$,
\begin{eqnarray} V_j^{-1}\equiv \{x\in {\bf R}^n;
x=p_j(s)+t{\mathbf v}_j, s\in S(0,k_\om),t\in (-d,0)\}\subset\om ,\\
V_j^{1}\equiv \{x\in {\bf R}^n; x=p_j(s)+t{\mathbf v}_j, s\in
S(0,k_\om),t\in (0,d)\}\subset {\bf R}^n\backslash\overline\om,
\end{eqnarray}
and
\begin{eqnarray} \{x\in {\bf R}^n; x=p_j(s)-d{\mathbf v}_j, s\in
S(0,k_\om)\}\subset \om. \end{eqnarray} We also notice  that
\begin{eqnarray}
V_j^{0}\equiv \{x\in {\bf R}^n; x=p_j(s), s\in S(0,k_\om)\}\subset
\p\om,
\end{eqnarray}

Now let $D$ and $B$ be bounded open subsets in ${\bf R}^n$ with
$B\subset\subset D$, and let $a,r$ and $k$ be positive constants
with $r>k$. Define ${\mathcal O}(a,r,k)$ to be {\it the family of
all open subsets contained in $B$ such that every $\om \in{\cal
O}(a,r,k)$ is of class $\cal C$ with $k_\om\geq k > 0, a_\om > 0,
r_\om\leq r$ ( where $k_\om, a_\om \ \mbox{and}\ r_\om$ are given
as in (2.1), (2.2)-(2.3) and (2.4) ). In what follows, we will
write ${\cal O}$ for ${\cal O}(a,r,k)$, for simplicity of
notation. Moreover, we require  that the family of functions
${\cal F}={\bigcup_{\om\in\cal O}} {\cal F}_\om$ is
equi-continuous and equi-bounded on $\overline {S(0,k)}$.} We
define the usual topology on $\cal O $ by the Hausdorff-Pompeiu
distance between the complementary sets (which are closed):
\begin{eqnarray}
\rho(\om_1,\om_2)=\mbox{dist}(\overline{B}\backslash
\om_1,\overline{B}\backslash \om_2)\ \ \mbox{for all} \ \
\om_1,\om_2\in\cal O. \end{eqnarray} Denote by Hlim, the limit in
terms of (2.11). The following lemmas will be used to study the
existence of our shape optimization problem in section 5.

\vskip 3mm

{\bf Lemma 2.1 }($\cite{kn:[9]}$)  Let
$\{\om_m\}_{m=1}^\infty\subset\cal O$. There exists a subsequence
$\{\om_{m_k}\}_{k=1}^\infty$ of $\{\om_{m}\}_{m=1}^\infty$ such
that
$$ \mbox{H}\!\!\!\lim_{k\rightarrow\infty}\om_{m_k}=\om_0\in\cal O.$$

\vskip 3mm

{\bf Lemma 2.2 }($\cite{kn:[6]}$)  ($\Gamma$-property for $\cal
O$) Let $\{\om_m\}_{m=1}^\infty\subset \cal O$ and let
$\om_0\in\cal O$ be such that
$\mbox{Hlim}_{m\rightarrow\infty}\om_m=\om_0$. Then for any open
subset K with $K\subset\subset \om_0$ (i.e.,$\overline K\subset
\om_0$), there is a natural number $m(K)$ such that for all $m\geq
m(K)$, $\overline K\subset \om_m$.

\vskip 3mm

{\bf Lemma 2.3 }($\cite{kn:[10]}$) ($ \hat{\Gamma}$-property for
$\cal O$) Let $\{\om_m\}_{m=1}^\infty\subset \cal O$ and
$\om_0\in\cal O$ be such that
$\mbox{Hlim}_{m\rightarrow\infty}\om_m=\om_0$. Then for any open
subset K with $\overline K\subset {\bf R}^n\backslash \om_0$,
there is a natural number $m(K)$ such that for any $m\geq m(K),
\overline K\subset D\setminus \om_m$.

\vskip 3mm

\def\theequation{3.\arabic{equation}}
\setcounter{equation}{0} \noindent {\large\bf3.\quad Divergence
free partition of vector-valued Sobolev functions ---an
application of the Hodge theorem}

\vskip 3mm

In this section, we shall apply the Hodge theorem to establish a
divergence free partition for vector-valued Sobolev functions in
${\bf R}^n, n\geq1$. We will refer the reader to  the book of
George de Rham $\cite{kn:[1111]}$  for a detailed accout on  the
Hodge Theory.

\vskip 3mm

We recall that for an open subset $\om\subset {\bf R}^n$,
 $D_0^1(\om)\equiv\{{\mathbf v}\in H_0^1(\om)^n;
\mbox{div}{\mathbf v}=0\}$.

\vskip 3mm

{\bf Theorem 3.1 } Let D be a bounded  open subset in ${\bf R}^n$
and  let ${\mathbf u}\in D_0^1(D)$.  Suppose that $U_i,
i=1,\cdots,m$, are bounded open subsets in ${\bf R}^n$ such that
$\bigcup_{j=1}^mU_j\subset\subset D$. Then for any open subset
$\hat U \subset\subset \bigcup_{j=1}^mU_j$, there are functions
$\{{\mathbf u}_j \}_{j=1}^m$ such that $${\mathbf u}=\sum_{j=1}^m
{\mathbf u}_j\ \ \mbox{ over }\ \ \hat U$$ and for  $j=1,\cdots,m,
\ {\mathbf u}_j\in D_0^1(U_j)$ and
 $$ \| {\mathbf u }_j \|_{H^1(U_j)}\leq \mbox{C} \|{\mathbf u}\|_{H^1(D)}, $$
where C is independent of ${\mathbf u} \in D_0^1(D)$.

\vskip 3mm

{\bf Proof of Theorem 3.1}\quad Write ${\mathbf
u}=(u_1,\cdots,u_n) $ and
\begin{eqnarray}\label{3.1} \omega =
\sum_{i=1}^n (-1)^{i+1} u_idx_1\wedge \cdots \wedge
\hat{dx_i}\wedge \cdots \wedge dx_n.\end{eqnarray}
Then
$d\omega=\mbox{div}{\mathbf u}\ dx_1\wedge\cdots\wedge dx_n=0$.
Hence, $\omega$ is a closed $(n-1)$-form over ${\bf R}^n$ with
$H^1({\bf R}^n)$-regular coefficients and $spt(\omega)\subset\overline D$,
where $spt(\omega)$ denotes the support of $\omega$.

Let $S^n$ be the standard sphere in ${\bf R}^{n+1}$ and let $g$ be
the standard metric on $S^n$. We fix a Sobolev space
$H^\alpha(S^n)$ of order $\alpha$ for  $\alpha \in {\bf R}$ by
using the metric $g$ and a certain partition of unity over $S^n$.
 Let $\pi: S^n\backslash N\rightarrow {\bf R}^n$ with
$N=(0,\cdots,0,1)$ be the stereographic projection and let

\begin{eqnarray*}
\omega^*=\begin{cases}{\pi^*(\omega),\ &$\mbox{for}\  p\in
S^n\backslash N$,\cr 0,\ &$\mbox{for}\ p=N.$}\end{cases}
\end{eqnarray*}

Since $\omega$ has compact support in ${\bf R}^n$, $\omega^*$ is an
$(n-1)$-form over $S^n$ with $H^1(S^n)$-regular coefficients. Moreover $\omega^*$ is
identically  zero in a small neighborhood of $N$. Notice that
$d\omega^*=0$.

Let $\Delta_{g,n-1}$ be the standard Betrami-Laplacian operator
over $S^n$ acting on the space of $(n-1)$-forms. Consider the
equation

\begin{eqnarray}
\Delta_{g,n-1}\sigma^*=\omega^*.
\end{eqnarray}
By the Hodge theorem, $\mbox{Ker}\Delta_{g,n-1}\cong H^{n-1}_{deRham
}(S^n,R)=\{0\}$. Hence, equation (3.2) admits a unique solution
$\sigma^*$ with the estimate
\begin{eqnarray*}
\|\sigma^*\|_{H^3(S^n)}\leq \mbox{C} \|\omega^*\|_{H^1(S^n)},
\end{eqnarray*}
where C is independent of the choice of $\omega^*$.

Notice that
$\Delta_{g,n-1}\sigma^*=dd^*\sigma^*+d^*d\sigma^*=\omega^*$ and
$d\omega^*=0$. Here $d^*$ is the Hilbert space adjoint of the
exterior differential operator. It follows that
\begin{eqnarray*}
(dd^*\sigma^*+d^*d\sigma^*,d^*d\sigma^*)=(\omega^*,d^*d\sigma^*)
=(d\omega^*,d\sigma^*)=0.
\end{eqnarray*}

Since $(dd^*\sigma^*,d^*d\sigma^*)=(d^2d^*\sigma^*,d\sigma^*)=0$,
we get $d^*d\sigma^*=0$. Thus $dd^*\sigma^*=\omega^*$. Write
$\tau^*=d^*\sigma^*$. Then $d\tau^*=\omega^*$ and
\begin{eqnarray*}
\|\tau^*\|_{H^2(S^n)}\leq
\mbox{C}\|\omega^*\|_{H^1(S^n)}\leq\mbox{C}\|\omega\|_{H^1(D)}.
\end{eqnarray*}
Here and in what follows, $C$ denotes a constant  independent of
the choice of $\omega$ with $spt\omega\subset \overline D$.
However the constant $C$ may be different in different context.

Let $\tau=(\pi^{-1})^*(\tau^*|_{S^n\backslash N})$. Then
\begin{eqnarray*}
d\tau=d((\pi^{-1})^*(\tau^*|_{S^n\backslash
N}))=(\pi^{-1})^*d(\tau^*|_{S^n\backslash
N})=(\pi^{-1})^*(\omega^*|_{S^n\setminus N})=\omega,
\end{eqnarray*}

and
\begin{eqnarray}
\|\tau\|_{H^2(D)}\leq\mbox{C}\|\tau^*\|_{H^2(S^n)}\leq\mbox{C}\|\omega\|_{H^1(D)}.
\end{eqnarray}

We next let $\{\chi_j\}_{j=1}^m$ be such that $\chi_j\in
C_0^\infty(U_j),0\leq\chi_j\leq 1$, and $\sum^m_{j=1}\chi_j=1$ over
$\hat U\subset\subset \bigcup^m_{j=1}U_j$. Set
$\tau_j=\chi_j\tau$. Then $spt\tau_j\subset U_j$ and
$$\sum^m_{j=1}d\tau_j=\omega \ \ \mbox{over} \ \ \hat U.$$
Moreover, $d\tau_j$ is an $(n-1)$-form over ${\bf R}^n$ with
$H^1$-regular coefficients and $spt (d\tau_j)\subset U_j$.
Write for each $j=1,\cdots,m$,
\begin{eqnarray}
d\tau_j=\sum^n_{i=1}(-1)^{i+1}b_{i,j}dx_1\wedge\cdots\wedge
\hat{dx_i}\wedge\cdots\wedge dx_n,
\end{eqnarray}
where $b_{i,j}\in H^1_0(U_j), \ i=1,\cdots,n.$

Let ${\mathbf u}_j=(b_{1,j},\cdots,b_{n,j})$. Then ${\mathbf
u}_j\in H^1_0(U_j)^n$. Since $0=d^2\tau_j=(\mbox{div}{\mathbf
u}_j)dx_1\wedge\cdots\wedge dx_n$, it follows that
$\mbox{div}{\mathbf u}_j=0$. Moreover, by (3.3) we have
\begin{eqnarray*}
\|{\mathbf u}_j\|_{H^1(U_j)}\leq\mbox{C}\|\tau_j\|_{H^2(U_j)}
\leq\mbox{C}\|\tau\|_{H^2(D)}\leq\mbox{C}\|\omega\|_{H^1(D)}
 \leq\mbox{C}\|{\mathbf u}\|_{H^1(D)}.
 \end{eqnarray*}

On the other hand, since
\begin{eqnarray*}
\sum^m_{j=1}d\tau_j=\omega \ \ \mbox{over}\ \ \hat U
\end{eqnarray*}
it follows from (3.1) and (3.4) that
\begin{eqnarray*}
&&\sum^m_{j=1}\sum^n_{i=1}(-1)^{i+1}b_{i,j}dx_1\wedge\cdots\wedge\hat{dx_i}\wedge\cdots\wedge
dx_n\\
&&=\sum^n_{j=1}(-1)^{i+1}u_idx_1\wedge\cdots\wedge\hat{dx_i}\wedge\cdots\wedge
dx_n\quad \quad \mbox{over}\quad \ \hat U, \end{eqnarray*} which
implies
\begin{eqnarray*}
u_i=\sum^m_{j=1}b_{i,j}\quad \ \mbox{over} \quad \hat U,\ \
i=1,\cdots,n.
\end{eqnarray*}

Thus, ${\mathbf u}=\sum_{i=1}^m{\mathbf u}_j$ over $\hat U$.  This
completed the proof.
$\endpf$

\vskip 3mm

{\bf Corollary 3.2 } Let $D\subset {\bf R}^n$ be a bounded open
subset and let ${\mathbf u}\in D^1_0(D)$. Suppose that
$\{U_j\}_{j=1}^m$ is an open covering of $\o{D}$ with
$\bigcup_{j=1}^m{U_j}\subset\subset B_{R}(0)$, the ball in ${\bf
R}^n$ centered at $0$ and of radius $R$.
Then there are functions $\{{\mathbf u}_j\}^m_{j=1}$     , such
that $${\mathbf u}=\sum^m_{j=1}{\mathbf u}_j  \ \mbox{over} \  D$$
and for all $j=1,\cdots,m$, ${\mathbf u}_j\in D^1_0(U_j) \
\mbox{and}\ \|{\mathbf u}_j\|_{H^1(U_j)}\leq\mbox{C}\|{\mathbf
u}\|_{H^1(D)}$. Here C is independent of ${\mathbf u}\in
D^1_0(D)$, but depends on $R$.

\vskip 3mm

{\bf Proof of Corollary 3.2}\quad Let $\omega,\pi, S^n, N,
\omega^*,\sigma^*,\tau^*$ and $\tau$ be given as in the proof of
theorem 3.1. We  then have $d\tau=\omega$ and
$\|\tau^*\|_{H^2(S^n)}\leq\mbox C\|\omega\|_{H^1(D)}$.
Let  $\chi_{B_R}\in C^\infty_0(B_R)$ be a cutting function with
$\chi_{B_R}=1$ over $\bigcup^m_{j=1}U_j$. Then
$$d(\chi_{B_R}\tau)=\omega \quad \ \mbox{over} \quad
\bigcup^m_{j=1}U_j$$ and
\begin{eqnarray}
\|\chi_{B_R}\tau\|_{H^2({\bf R}^n)}\!\!\!\!
 && \leq
\mbox{C}\|\pi^*(\chi_{B_R}\tau)\|_{H^2(S^n)} \leq
\mbox{C}\|\tau^*\|_{H^2(S^n)} \cr
 &&\leq\mbox{C}\|\omega^*\|_{H^1(S^n)} \leq \mbox{C} \| \omega
 \|_{H^1(D)}.
\end{eqnarray}

Now let $\{\chi_j\}_{j=1}^m$ be such that $0\leq \chi_j\leq
1,\chi_j\in C_0^\infty(U_j),\sum_{j=1}^m\chi_j=1$ over $D$. Then
$$d(\chi_{B_R}\tau\sum_{j=1}^m\chi_j)=\omega\quad \ \mbox{over}\quad D.$$
Let $\tau_j=\chi_{B_R}\chi_j\tau$. Then $spt\tau_j\subset U_j$.
Thus it follows from (3.5) that $d\tau_j$ is an $(n-1)$-form with
$H^1({\bf R}^n)$-smooth coefficients and
$$\|d\tau_j\|_{H^1({\bf R}^n)}\leq\mbox{C}\|\omega\|_{H^1(D)}.$$
Write
$$d\tau_j=\sum_{i=1}^n(-1)^{i+1}b_{i,j}dx_1\wedge\cdots\wedge
\hat{dx_i}\wedge\cdots\wedge dx_n,j=1,\cdots,m,$$ and let
${\mathbf u}_j=(b_{1,j},\cdots,b_{n,j})$ for $j=1,\cdots,m$. Then
by the same argument as that in the proof of Theorem 3.1, we get
\begin{eqnarray*}&&{\mathbf u}=\sum_{j=1}^m{\mathbf u}_j\quad \
\mbox{over}\quad D,\\
&& {\mathbf u}_j\in D_0^1(U_j),j=1,\cdots,m,
\end{eqnarray*}
and
$$\|{\mathbf u}_j\|_{H^1(U_j)}\leq\mbox{C}\|{\mathbf
u}\|_{H^1(D)},j=1,\cdots,m.$$ This completes the proof. $\endpf$
\bigskip

\def\theequation{4.\arabic{equation}}
\setcounter{equation}{0} \noindent {\large\bf4.\quad A localization
property  of the divergence free decomposition}

\vskip 3mm

In this section, we shall prove a localization property of the
divergence free decomposition established in section 3. Based on
it, we obtain some identities for certain fundamental function
spaces related to the Navier-Stokes equations.

\vskip 3mm

{\bf Theorem 4.1 } Let $D$ be a bounded open subset in ${\bf
R}^2$, and let $\Omega=\bigcup_{j=1}^N\Omega_j\subset\subset D$ be
a finite union of open and connected subsets in ${\bf R}^2$ with
$\overline{\Omega_i}\bigcap\overline{\Omega_j}=\emptyset$ for
$i\not=j$. Assume that $\{U_j\}_{j=1}^m$ be an open covering of
$\p\Omega$ in ${\bf R}^2$ with $\bigcup_{j=1}^mU_j\subset\subset
D$. Let ${\mathbf u}\in D_0^1(D) $ with ${\mathbf u}=0$ over
$\Omega$. Then there are functions $\{{\mathbf u}_j\}_{j=1}^m$
such that
$${\mathbf u}=\sum_{j=1}^m{\mathbf u}_j$$
over a certain neighborhood $\hat U$ of $\p\Omega$ in ${\bf R}^2$
with $\hat U\subset\subset \bigcup_{j=1}^m U_j$, and for all
$j=1,\cdots,m $, ${\mathbf u}_j\in D_0^1(U_j), {\mathbf u}_j=0$
over $\Omega$ and
$$\|{\mathbf u}_j\|_{H^1(U_j)}\leq \mbox{C}\|{\mathbf
u}\|_{H^1(D)},$$ where C is a positive constant independent of
${\mathbf u}\in D_0^1(D)$.

 \vskip 3mm

{\bf Proof of Theorem 4.1} \quad Step 1. We shall first prove the theorem in the
case that $\om$ is connected.

Let $\omega,\omega^*,\tau^*$ and $\pi$ be given as in the proof of
Theorem 3.1. We  then have $\omega^*=d\tau^*$ and
\begin{eqnarray}
\|\tau^*\|_{H^2(S^2)}\leq\mbox{C}\|\omega\|_{H^1(S^2)}\leq\mbox{C}\|\omega\|_{H^1(D)}.
\end{eqnarray}

Since ${\mathbf u}=0$ over $\om$, it follows that $\omega=0$ over
$\om$, which implies  that $\omega^*=0$ over $\pi^{-1}(\om)$.
Hence, $d\tau^*=0$ over $\pi^{-1}(\om)$. Because $\tau^*$ is a
0-form and $\om$ is connected, there is a constant $c_0\in {\bf
R}$ such that
\begin{eqnarray}
\tau^*=c_0\quad \ \mbox{over}\quad \pi^{-1}(\om).
\end{eqnarray}
By (4.1) and (4.2), we obtain
\begin{eqnarray}
|c_0|\leq\mbox{C}\|\omega\|_{H^1(D)}.
\end{eqnarray}

Let $\tilde{\tau}^*=\tau^*-c_0$. Then $\tilde{\tau}^*=0$ over
$\pi^{-1}(\om)$ and
$d\tilde{\tau}^*=d(\tau^*-c_0)=d\tau^*=\omega^*$, moreover, it
follows from (4.1) and (4.3) that
$$\|\tilde{\tau}^*\|_{H^2(S^2)}\leq C \|\omega\|_{H^1(D)}.$$
Let $\tilde{\tau}=(\pi^{-1})^*(\tilde{\tau}^*|_{S^2\setminus N})$.
Then $\tilde{\tau}=0$ over $\om$, $d\tilde{\tau }=\omega$ and
\begin{eqnarray}
\|\tilde{\tau}\|_{H^2(D)}\leq \mbox{C}
\|\pi^*\tilde{\tau}\|_{H^2(S^2)}\leq\mbox{C}\|\tilde{\tau}^*\|_{H^2(S^2)}\leq
\mbox{C}\|\omega\|_{H^1(D)}.
\end{eqnarray}

Let $\hat U$ be any neighborhood of $\p\om$ such that $\hat U
\subset\subset \bigcup^m_{j=1}U_j$ and let $\{\chi_j\}^m_{j=1}$ be
such that $0\leq\chi_j\leq 1,\chi_j\in C_0^1(U_j)$  and
$\sum^m_{j=1}\chi_j=1$ over $\hat U$.  Then we have
$$d(\sum^m_{j=1}\chi_j\tilde{\tau})=\omega\quad\ \mbox{over} \quad
\hat U.$$ Let $\tilde\tau_j=\chi_j\tilde \tau$,  $j=1,\cdots,m$.
Then $\tilde\tau_j\in H_0^1(U_j)$, $\tilde\tau_j=0$ over $\om$ and
$$\sum_{j=1}^m d\tilde\tau_j =\omega  \ \mbox{over}\ \hat U.$$

Write, for $j=1,\cdots,m, d\tilde\tau_j=b_{1,j}dx_2-b_{2,j}dx_1 $
for certain $b_{1,j},b_{2,j}\in H_0^1(U_j)$ with $b_{i,j}=0$ over
$\om$. Let ${\mathbf u}_j=(b_{1,j},b_{2,j}),j=1,\cdots,m$. Then
one can easily check that for all $j=1,\cdots,m,{\mathbf u}_j=0$
over $\om$, ${\mathbf u}_j\in D_0^1(U_j)$ and ${\mathbf
u}=\sum_{j=1}^m{\mathbf u}_j$ over $\hat U$.
  Moreover, it follows from (4.4) that
\begin{eqnarray*}
\|{\mathbf u}_j\|_{H^1(U_j)} &&
\leq\mbox{C}(\|b_{1,j}\|_{H^1(U_j)}+\|b_{2,j}\|_{H^1(U_j)})\\
&&\leq\mbox{C}\|d\tilde{\tau_j}\|_{H^1(U_j)}\leq\mbox{C}\|\tilde{\tau}\|_{H^2(D)}\\
&&\leq\mbox{C}\|\omega\|_{H^1(D)}\leq\mbox{C}\|{\mathbf
u}\|_{H^1(D)}.
\end{eqnarray*}
This proves the theorem in case  $\om$ is connected.

\vskip 3mm

Step 2. We next prove the theorem in case $\om$ is
multiple-connected. For brevity of notation, we assume that $\om$
is 2-connected. (The general case can be similarly done through an
induction argument.)  Then $\om=\om_1\bigcup\om_2$ with
$\overline{\om}_1\bigcap\overline{\om}_2=\emptyset$, $\om_1$ and
$\om_2$ being connected.

By refining the covering $\{U_j\}^m_{j=1}$ of $\p\om$, we may
assume, without loss of generality, that there are open subsets
$\{U_{j,1}\}^{m}_{j=1}$ and $\{U_{j,2}\}^{m}_{j=1}$
 such that
\begin{eqnarray*}
&&\bigcup_{j=1}^{m}U_{j,1}\supset\p\om_1,\
\bigcup_{j=1}^{m}U_{j,2}\supset \p\om_2,\\
&&U_{j,i}\subset U_j\ \ \mbox{for}\ \ j=1,\cdots,m,\ i=1,2,\\
&&U_{l,1}\bigcap U_{j,2}=\emptyset \ \ \mbox{for }\ \
l,j=1,\cdots,m.
\end{eqnarray*}
Let $\hat U_1$ and $\hat U_2$ be neighborhoods of $\p\om_1$ and
$\p\om_2$, respectively, such that
$$\hat U_1\subset\subset \bigcup_{j=1}^{m_1}U_{j,1},\ \
\hat U_2\subset\subset \bigcup_{j=1}^{m_2}U_{j,2}.$$ Since
${\mathbf u}\in D_0^1(D)$ and ${\mathbf u}=0$ over $\om_1$, it
follows from Step 1 that there are functions $\{{{\mathbf
u}_{j,1}}\}^{m}_{j=1}$ such that  ${\mathbf
u}=\sum^{m}_{j=1}{\mathbf u}_{j,1}$ over $\hat U_1$ and for each
$j=1,\cdots,m$, ${\mathbf u}_{j,1}\in D_0^1(U_{j,1}),{\mathbf
u}_{j,1}=0$ over $\om_1$, and
$$\|{\mathbf u}_{j,1}\|_{H^1(U_{j,1})}\leq \mbox{C}\|{\mathbf
u}\|_{H^1(D)}.$$ Notice that for  $j=1,\cdots,m,\ {\mathbf
u}_{j,1}=0$ over $\om_2$, for $\om_2\bigcap U_{j,1}=\emptyset$ and
${\mathbf u}_{j,1}\in D_0^1(U_{j,1})$.

Now let ${\mathbf v}={\mathbf u}-\sum^{m}_{j=1}{\mathbf u}_j$.
Then ${\mathbf v}\in D_0^1(D)$ and ${\mathbf v}=0$ over $\om_2$.
Again, by Step 1, there are functions $\{{\mathbf
u}_{j,2}\}_{j=1}^{m}$ such that for each $j=1,\cdots,m$, ${\mathbf
u}_{j,2}\in D_0^1(U_{j,2}),{\mathbf u}_{j,2}=0$ over $\om_2$,
$$\|{\mathbf u}_{j,2}\|_{H^1(U_{j,2})}\leq\mbox C\|{\mathbf
v}\|_{H^1(D)}\leq\mbox C\|{\mathbf u}\|_{H^1(D)},$$ and ${\mathbf
v}=\sum_{l=1}^{m}{\mathbf u}_{j,2}$ over $\hat U_2$.  Notice also
that for  $j=1,\cdots,m,\ {\mathbf u}_{j,2}=0$ over
 $\om_1$, for $U_{j,2}\bigcap\om_1=\emptyset$ and ${\mathbf u}_{j,2}\in
D_0^1(U_{j,2})$.

Finally, let ${\bf u}_l={\bf u}_{l,1}+{\bf u}_{l,2},\
l=1,\cdots,m$.  Then it is clear that
$${\mathbf u}=\sum_{l=1}^m{\mathbf u}_l\ \ \mbox{ over }\ \ \hat U\equiv\hat U_1\bigcup \hat U_2,$$
and for all $l=1,\cdots,m,\ {\mathbf u}_l\in D_0^1(U_l),\ {\mathbf
u}_l=0$ over $\om$ and $\|{\mathbf u}_l\|_{H^1(U_l)}\leq
\mbox{C}\|{\mathbf u}\|_{H^1(D)}$.  This completes the proof.
$\endpf$

\vskip 3mm

By the almost same argument as above with some slight
modification, we can get the following.

\vskip 3mm

{\bf Theorem 4.1' }  Let $D$ be a bounded open subset in ${\bf
R}^2$, and let $\om=\bigcup_{j=1}^N\om_j\subset\subset D$ be a
finite union of open and connected subsets with
$\overline{\om_i}\bigcap\overline{\om_j}=\emptyset$ for $i\not=j$.
Suppose that $\Gamma$ is a compact subset of $\p\om$ and
$\{U_j\}_{j=1}^m $ is an open covering of $\Gamma$ in ${\bf R}^2$
with $\bigcup_{j=1}^mU_j\subset\subset D$. Let ${\mathbf u}\in
D_0^1(D)$ with ${\mathbf u}=0$ over $\om$. Then there are
functions $\{{\mathbf u}_j\}_{j=1}^m$ such that
$${\mathbf u}=\sum_{j=1}^m{\mathbf u}_j$$
over a certain neighborhood $\hat U$ of $\Gamma$ in ${\bf R}^2$
with $\hat U\subset\subset \bigcup_{j=1}^m U_j$, and for all
$j=1,\cdots,m, {\mathbf u}_j\in D_0^1(U_j), {\mathbf u}_j=0$ over
$\om$ and
$$\|{\mathbf u}_j\|_{H^1(U_j)}\leq\mbox{C}\|{\mathbf
u}\|_{H^1(D)},$$ where C is a positive constant independent of
${\mathbf u}\in D_0^1(D)$.

\vskip 3mm

We next prove a version of  Theorem 4.1 for domains in ${\bf
R}^3$. To this aim, we need to recall some results from Algebraic
Topology. Two good references for the topological results which we
will quote here are $\cite{kn:[11]}$ and $\cite{kn:[111]}$.

Suppose that $\om$ is a bounded open subset in ${\bf R}^n,n\geq
1$, with $\p\om$ a compact topological manifold of dimensional
$n-1$. By making use of results in Algebraic Topology( see pp.227,
(26.17.8), $\cite{kn:[11]}$), there is a bounded domain
$\om'\subset\subset {\bf R}^n$ with smooth boundary $\p\om'$ such
that $\om'\supset\supset\om$ and $\om'$ can be retracted to $\om$.
Namely, there is a continuous map $\gamma: \om'\rightarrow
\overline\om$ with $\gamma|_{\overline\om}=id$ ( the identity map
). This implies that the homological group $H_1(\om, {\bf R})$ is
of finite dimension and can be naturally embedded as a vector
space of $H_1(\om',{\bf R})$. Let
$\{[\gamma_1],\cdots,[\gamma_k]\}$ be a basis of $H_1(\om,{\bf
R})$ with $\gamma_j: [0,1]\rightarrow\om$ a smooth  closed Jordan
curve for each j (see pp.63, (12.1), $\cite{kn:[11]}$). We can
then extend $\{[\gamma_1],\cdots,[\gamma_k]\}$ to a basis
$\{[\gamma_1],\cdots,[\gamma_k],\cdots,[\gamma_M]\}$ of
$H_1(\om',{\bf R})$. Let $\{[\tau_1],\cdots,[\tau_M]\}$ be the
dual basis of $\{[\gamma_1],\cdots,[\gamma_k],\cdots,[\gamma_M]\}$
in $H^1_{de Rham}(\om',{\bf R})$. Then
$\{[\tau_1],\cdots,[\tau_k]\}$, when restricted to $\om$, gives a
dual basis of $\{[\gamma_1],\cdots,[\gamma_k]\}$, namely,
$\displaystyle \int_{\gamma_l}\tau_j=\delta_j^l$. Here
$\{\tau_j\}$ are smooth 1-forms over $\om'$.

\vskip 3mm

{\bf Lemma 4.2 } Under the above notations, further assume that
$\om\subset {\bf R}^3$ is connected and   of class $\cal C$. Let
$\tau$ be a 1-form over $\om'$ with $H^2(\om')$-regular
coefficients. Also assume that $d\tau=0$ over $\om$. Then there
exists $h\in H_0^2(\om')$ and $\{c_j\}_{j=1}^k\subset {\bf R}$
such that
\begin{eqnarray*}
 \tau=dh+\sum_{j=1}^kc_j\tau_j \quad \ \mbox{over} \quad \om,
\end{eqnarray*}
 and
\begin{eqnarray*}
  \|h\|_{H^2(\om')}\leq
\mbox{C}\|\tau\|_{H^2(\om')},\  |c_j |\leq \mbox
C\|\tau\|_{H^2(\om')},\ j=1,\cdots,k.
\end{eqnarray*}
Here C is a positive constant depending only on $\om$ and $\om'$.

\vskip 3mm

{\bf Proof of Lemma 4.2}\quad Let $\psi_0\in C_0^\infty(-1,1) $ be
such that $\psi_0\geq0$ and $\int_{{\bf R}^3}\psi(x)dx=1$ where
$\psi(x)\equiv \psi_0(|x|^2)$, and let $\tau_\delta(x)=\int_{{\bf
R}^3}\tau(x-y)\psi_\delta(y)dy$, where $\delta>0$ and
$\psi_\delta(y)=\delta^{-3}\psi(\delta^{-1}y)$. Then by the
Friedrich smoothing lemma, for any $\delta_0<<1$,
$\tau_\delta\rightarrow\tau$ in $H^2(\om'_{\delta_0})$ as
$\delta\rightarrow 0$ where
$\om'_{\delta_0}\equiv\{x\in\om';\mbox{dist}(x,\p\om')>\delta_0\}$.
Moreover, $\tau_\delta\in\Lambda^1(\om'_\delta)$,  i.e.,
$\tau_\delta$ is a smooth 1-form over $\om'_\delta$ and
$d\tau_\delta=0$ over $\om_\delta$, where
$\om_\delta\equiv\{x\in\om; \mbox{dist}(x,\p\om)>\delta\}$.

Next, since $\tau$ has $H^2(\om')$-regular coefficients, for each
piecewise smooth curve $\gamma$ in $\om'$, the restriction of
$\tau$ on $\gamma$ is well-defined and is $L^2$-integrable by the
trace theorem. Hence, we can well define  $\int_\gamma\tau$,
moreover, by the trace theorem,
$$|\int_\gamma\tau|\leq \mbox{C}(\gamma)\|\tau\|_{H^1(\om')},$$
where C$(\gamma)$ depends on $\gamma$.

Now let $c_j=\int_{\gamma_j}\tau,\ \ j=1,\cdots,k$, and   $\hat
\tau=\tau-\sum^k_{j=1}c_j\tau_j$.  It is clear that
$$|c_j|\leq\mbox{C}(\gamma_j)\|\tau\|_{H^1(\om')}.$$

We claim that for each piecewise smooth curve $\gamma:
[0,1]\rightarrow\om$, $\int_\gamma\hat\tau$ depends only on
$\gamma(0)$ and $\gamma(1)$.

Indeed, for each $\gamma': [0,1]\rightarrow\om$ piecewise smooth
curve with $\gamma'(0)=\gamma(0), \gamma'(1)=\gamma(1)$, since
$\gamma'-\gamma$ is homotopical to $\sum_{j=1}^kd_j\gamma_j$ for
a certain choice of $\{d_j\}_{j=1}^k\subset {\bf R}$, we get
\begin{eqnarray*}
\int_{\gamma'-\gamma}(\tau_\delta-\sum_{i=1}^kc_i\tau_i)
 &&=\sum_{j=1}^k\int_{\gamma_j}d_j(\tau_\delta-\sum_{i=1}^kc_i\tau_i)\\
 &&=\sum_{i,j=1}^k(d_j\int_{\gamma_j}\tau_\delta-c_id_j\delta_i^j)\\
 &&=\sum_{j=1}^kd_j(\int_{\gamma_j}\tau_\delta-c_j).
 \end{eqnarray*}

Letting $\delta\rightarrow0$, we get
$$\int_{\gamma'-\gamma}(\tau-\sum_{i=1}^kc_i\tau_i)=\sum_{j=1}^kd_j(c_j-c_j)=0.$$
Now we fix $p_0\in\om$. Since $\om$ is connected, for any
$p\in\om$, there is a piecewise smooth curve $\gamma:
[0,1]\rightarrow\om$ with $\gamma(0)=p_0,\gamma(1)=p$. We define
$h: \om\rightarrow {\bf R}$ by
$$h(p)=\int_{\gamma}(\tau-\sum_{j=1}^kc_j\tau_j)=\int_\gamma\hat\tau.$$
It is clear that $h$ is well-defined over $\om$. By the standard
arguments in calculus, one can easily verify that $dh=\hat\tau$
over $\om$ and thus $h\in H^3(\om)$. Moreover,
$$\|h\|_{H^3(\om)}\leq \mbox{C}\|\tau\|_{H^2(\om')}.$$
Through  the proof, as mentioned before, C denotes a positive
constant depending only on $\om$ and $\om'$, which may be different in different context.

Next, we are going to extend $h$ from $\om$ to $\om'$ so that the
extension is the desired function. To this end, we observe first
that since $\om$ is of class $\cal C$, there exist open sets
$\{V_j\}_{j=1}^m$ such that (2.5)-(2.10) are fulfilled. Moreover,
by taking $d$ small enough, we may have
$$\bigcup_{i=1}^mV_i\subset\subset\om'.$$
Now we define the extension of $h$ on $V_j$ as follows:
$$\tilde{h}_j(p_j(s)+t{\mathbf v}_j)=\int_{C_{p_j(s)+t{\mathbf
v}_j}}\hat\tau+h(p_j(s)-d{\mathbf v}_j),$$ for  $s\in
S(0,k_\om)$ and $t\in(-d,d)$, where $C_{p_j(s)+t{\mathbf v}_j}$
denotes the segment connecting $p_j(s)-d{\mathbf v}_j$ to
$p_j(s)+t{\mathbf v}_j$. It is clear that $\tilde h_j$ is
well-defined. By the trace theorem and a direct computation, we
have $\tilde h_j \in H^2(U_j)$ with $\|\tilde h_j
\|_{H^2(U_j)}\leq\mbox C\|\tau\|_{H^2(\om')}$. Moreover, by
(2.5)-(2.10), the whole segment $C_{p_j(s)+t{\mathbf v}_j}$ lies
in $\om$ if $p_j(s)+t{\mathbf v}_j\in\om\bigcap V_j$. Thus for
$p_j(s)+t{\mathbf v}_j\in\om\bigcap V_j$,
\begin{eqnarray*}
\tilde h_j(p_j(s)+t{\mathbf v}_j)=\int_{C_{p_j(s)+t{\mathbf
v}_j}}dh+h(p_j(s)-d{\mathbf v}_j)=h(p_j(s)+t{\mathbf v}_j),
\end{eqnarray*}
i.e., $\tilde h_j=h$ over $\om\bigcap V_j$.

Let $\hat U$ be a neighborhood of $\p\om$ such that $\hat U
\subset\subset\bigcup_{j=1}^mV_j$ and $\{\chi_j\}_{j=1}^m$ be such
that $\chi_j \in C_0^\infty(V_j),0\leq\chi_j\leq 1$ for
$j=1,\cdots,m$, and $\sum_{j=1}^m\chi_j=1$ over $\hat U$. Set
\begin{eqnarray*}
\tilde{\tilde h}_j(x)=
\begin{cases}{\tilde h_j(x),\ &  $ x\in V_j$,\cr
 0\ , &  $x\in \om'\setminus\overline V_j$.}
 \end{cases}
\end{eqnarray*}
Then $\chi_j\tilde{\tilde h}_j \in H^2(\om')$ and
$$\|\chi_j\tilde{\tilde h}_j\|_{H^2(\om')}\leq\mbox
C\|\tau\|_{H^2(\om')}$$ Let $\tilde
h=\sum_{j=1}^m\chi_j\tilde{\tilde h}_j$, then $\tilde h(x)=h(x)$
over $\om\bigcap\hat U$ and
$$\|\tilde h\|_{H^2(\om')}\leq \mbox C\|\tau\|_{H^2(\om')}$$
Finally, we set
\begin{eqnarray*}
\hat h(x)=
\begin{cases}{h(x),& $x\in \om$,\cr \tilde h(x),&$x\in \hat U \bigcup
(\om'\setminus\overline\om)$.}
\end{cases}
\end{eqnarray*}
Then one can easily verify that $\hat h\in H^2(\om'), \hat h=h$
over $\om$ and
$$\|\hat h\|_{H^2(\om')}\leq \mbox C\|\tau\|_{H^2(\om')}.$$
This completes the proof of Lemma 4.2. $\endpf$

\vskip 3mm

{\bf Remark 4.3 } Let $D\subset {\bf R}^n$ be an open subset such
that $\om\subset\subset D$. Then $\om'$ in the above lemma can be
taken as $\om'\subset\subset D$.

\vskip 3mm

{\bf Theorem 4.4 } Let $D$ be a bounded open subset in ${\bf R}^3$
and let $\om\subset\subset D$ be  of class $\cal C$. Assume that
$\om'\supset\supset \om$ is given as in Lemma 4.2 such that
$\om'\subset\subset D$. Let $\{U_j\}_{j=1}^m$ be an open covering
of $\p\om$ in ${\bf R}^3$ with $\bigcup_{j=1}^m U_j \subset
\subset \om'$ and let ${\mathbf u}\in D_0^1(D)$ with ${\mathbf
u}=0$ over $\om$. Then there are functions $\{{\mathbf
u}_j\}_{j=1}^m$ such that
$${\mathbf u}=\sum_{j=1}^m{\mathbf u}_j$$
over a certain neighborhood $\hat U$ of $\p\om$ in ${\bf R}^3$
with $\hat U\subset\subset \bigcup_{j=1}^mU_j$, and for all
$j=1,\cdots,m$, ${\mathbf u}_j\in D_0^1(U_j)$, ${\mathbf u}_j=0$
over $\om$ and
$$\|{\mathbf u}_j\|_{H^1(U_j)}\leq \mbox{C}\|{\mathbf
u}\|_{H^1(D)},$$ where C is a positive constant independent of
${\mathbf u}\in D_0^1(D)$.

\vskip 3mm

{\bf Proof of Theorem 4.4}\quad Since $\om$ is of class $\cal C$,
we have $\om=\bigcup_{j=1}^N\om_j$ with $\om_i$ being open and
connected and $\overline{\om_i}\bigcap\overline{\om_j}=\emptyset$
for $i\not=j$.  We shall prove the theorem only in the case when
$\om$ is connected. The rest follows from the identical argument
as that in Step 2 of the proof of Theorem 4.1.

Let $\omega,\pi,\omega^*,\sigma^*,\tau^* $ and $\tau $ be given as
in the proof of   Theorem 3.1. We have $$d\tau=\omega \quad
\mbox{and}\quad \|\tau\|_{H^2(D)}\leq \mbox
C\|\omega\|_{H^1(D)}.$$ Here and throughout  the proof, $C$
denotes a  positive constant independent of $\omega$,
whenever $spt\omega\subset \overline D$. However $C$ may be different in different context.
Thus $\tau|_{\om'}$, the
restriction of $\tau$ over $\om'$, is an 1-form over $\om'$ with
$H^2(\om')$-regular coefficients.

Since ${\mathbf u}=0$ over $\om$, it follows that $d\tau=\omega=0$
over $\om$. Then by Lemma 4.2, there are $h\in H_0^2(\om')$,
$\{\tau_j\}_{j=1}^k\subset\Lambda^1(\om')$ and
$\{c_j\}_{j=1}^k\subset {\bf R}$ such that
\begin{eqnarray*}
&&\tau=dh+\sum_{j=1}^kc_j\tau_j \quad \quad \mbox{over} \quad
\om,\\
&&\|h\|_{H^2(\om')}+\sum_{j=1}^k|c_j|\leq\mbox
C\|\tau\|_{H^2(\om')}\leq\mbox C\|\tau\|_{H^2(D)}.
\end{eqnarray*}
Let $\hat U$ be a neighborhood of $\p\om$ in ${\bf R}^3$ such that $\hat
U\subset\subset \bigcup_{j=1}^mU_j$ and $\{\chi_j\}_{j=1}^m$ be
such that $\chi_j\in C_0^\infty(U_j), 0\leq\chi_j\leq 1$ and
$\sum_{j=1}^m\chi_j=1$ over $\hat U$.

Let
$$\tau^j=\chi_j(\tau-dh-\sum_{i=1}^kc_i\tau_i),\ \
j=1,\cdots,m.$$
 Then $\tau^j$ is a  1-form over ${\bf R}^3$ and
$$d\tau^j=d\chi_j\wedge(\tau-dh-\sum_{i=1}^kc_i\tau_i)+\chi_jd\tau.$$
It is  obvious that
$$\sum_{j=1}^md\tau^j=\omega\ \ \mbox{ over } \ \ \hat U,$$
and for all $j=1,\cdots,m,\  d\tau^j\in H_0^1(U_j),\ d\tau^j=0$
over $U_j\bigcap\om$. Moreover, it follows from the estimate in
Lemma 4.2 that for all $j=1,\cdots,m$,
$$\|d\tau^j\|_{H^1(U_j)}\leq
\mbox{C}\|\tau\|_{H^2(D)}\leq\mbox{C}\|\omega\|_{H^1(D)}.$$
 Now we write
$$d\tau^j=\sum_{i=1}^3(-1)^{i+1}b_{i,j}dx_1\wedge\cdots\wedge
\hat{dx_i}\wedge\cdots\wedge dx_3,\ j=1,\cdots,m,$$
 and let
${\mathbf u}_j=(b_{1,j},b_{2,j},b_{3,j}),\ j=1,\cdots,m$. Then one
can check easily that
$${\mathbf u}=\sum_{j=1}^m{\mathbf u}_j\ \ \mbox{ over }\ \ \hat
U,$$
 and for all $j=1,\cdots,m,\ {\mathbf u}_j\in D_0^1(U_j),\ {\mathbf
 u}_j =0 $ over $\om$ and
 $$\|{\mathbf u}_j\|_{H^1(U_j)}\leq\mbox{C}\|{\mathbf
 u}\|_{H^1(D)}.$$
  This completes the proof of Theorem 4.4.
 $\endpf$

 \vskip 3mm

By the same argument as above with some slight modification, we
can get the following.

\vskip 3mm

{\bf Theorem 4.4' } Let $D$ be a bounded open subset in ${\bf
R}^3$ and let $\om\subset\subset D$ be of class $\cal C$. Assume
that $\om'\supset\supset \om$ is given as in Lemma 4.2 such that
$\om'\subset\subset D$. Let $\Gamma$ be a compact subset of
$\p\om$ and let $\{U_j\}_{j=1}^m$ be an open covering of $\Gamma$
in ${\bf R}^3$ with $\bigcup_{j=1}^m U_j\subset\subset \om'$.
Suppose that ${\mathbf u}\in D_0^1(D)$ with ${\mathbf u}=0$ over
$\om$. Then there are functions $\{{\mathbf u}_j\}_{j=1}^m$ such
that
$${\mathbf u}=\sum_{j=1}^m{\mathbf u}_j$$
 over a certain neighborhood $\hat U$ of $\Gamma$ in ${\bf R}^3$
   with $\hat U\subset\subset \bigcup_{j=1}^mU_j$,  and
 for all $j=1,\cdots,m,\ {\mathbf u}_j\in D_0^1(U_j),\ {\mathbf
 u}_j=0$ over $\om$ and
$$\|{\mathbf u}_j\|_{H^1(U_j)}\leq\mbox{C}\|{\mathbf
u}\|_{H^1(D)},$$
 where C is a positive constant independent of ${\mathbf u}\in
 D_0^1(D)$.

 \vskip 3mm

 We next apply the above results to obtain some characterizations
 for some Sobolev function spaces.

\vskip 3mm

{\bf Theorem 4.5 } Let $D\subset {\bf R}^n, \ n=2,3$, be a bounded
open Lipschitz subset and let $\om$ be of class $\cal C$ such that
$\om\subset\subset D$. Then
$H_{0,\sigma}^1(D\setminus\overline\om)=\{{\mathbf v}\in
H_{0,\sigma}^1(D);\ {\mathbf v}=0\ \ \mbox{a.e. in }\  \om\}$

\vskip 3mm

{\bf Proof of Theorem 4.5}\quad It is clear that
$H_{0,\sigma}^1(D\setminus\overline\om)\subset\{{\mathbf v}\in
H_{0,\sigma}^1(D);\ {\mathbf v}=0\ \ \mbox{a.e. in }\ \ \om\}$.

Let ${\mathbf u}\in H_{0,\sigma}^1(D)$ satisfying ${\mathbf u}=0$
a.e. in $\om$. Let  $\om'$  be given as in Lemma 4.2 such that
$\om\subset\subset\om'\subset\subset D$. Let $\{V_j\}_{j=1}^m $ be
the covering of $\p\om$  as given by (2.5) and (2.6). By taking
$d$ in (2.6) sufficiently small, we have
$$ \bigcup_{j=1}^mV_j\subset\subset \om'.$$
Then by Theorem 4.1 and Theorem 4.4, there are functions
$\{{\mathbf u}_j\}_{j=1}^m$ and a neighborhood $\hat U$ of $\p\om$
in ${\bf R}^n$ such that
$${\mathbf u}=\sum_{j=1}^m{\mathbf u}_j\ \ \mbox{ over } \ \ \hat
U$$
 and for all $j=1,\cdots,m,\ {\mathbf u}_j\in D_0^1(V_j),\
 {\mathbf u}_j=0$ over $\om$.

Since $V_j\subset\subset D$ and ${\mathbf u}_j\in D_0^1(V_j)$, it holds that
 ${\mathbf u}_j\in D_0^1(D)$. Since $D$ is Lipschitz, we
have $D_0^1(D)=H^1_{0,\sigma}(D)$ (see Lemma 1.2.1, Chapter 3,
$\cite{kn:[7]}$). Hence, ${\mathbf u}_j\in H_{0,\sigma}^1(D)$ for
 $j=1,\cdots,m$.

Now let ${\mathbf u}_0={\mathbf u}-\sum_{j=1}^m{\mathbf u}_j$.
Then ${\mathbf u}_0\in H_{0,\sigma}^1(D),{\mathbf u}_0=0$ a.e. in
$\om\bigcup\hat U$. Since $\hat U\bigcup \om\supset\supset \om$,
there exists an open  Lipschitz subset $\hat V$ such that $\hat
U\bigcup \om \supset\supset\hat V\supset\supset\om$. Hence,
${\mathbf u}_0=0$ a.e. in $\hat V$ which implies by the trace
theorem that ${\mathbf u}_0\in H_0^1(D\setminus\overline{\hat
V})^n$. Since $D\setminus \overline{\hat V}$ is Lipschitz  and
div${\mathbf u}_0=0$, we have ${\mathbf u}_0\in
H_{0,\sigma}^1(D\setminus \overline{\hat V})$ (see Lemma 1.2.1,
Chapter 3, $\cite{kn:[7]}$). Thus ${\mathbf u}_0\in
H_{0,\sigma}^1(D\setminus\overline\om)$.

We next claim that ${\mathbf u}_j\in
H_{0,\sigma}^1(D\setminus\overline\om)$ for each  $j=1,\cdots,m$.
To this aim, we fix $j\in\{1,\cdots,m\}$ and define ${\mathbf
u}_j^t$ for $t\in(0,r)$, where
 $r=\f14\mbox{min}(d,
\mbox{dist}(spt{\mathbf u}_j,\p V_j)),$ by ${\mathbf
u}_j^t(x)={\mathbf u}_j(x-t{\mathbf v}_j)$. ( We recall that
${\mathbf v}_j$ is given by (2.6). ) One can easily check  that
div${\mathbf u}_j^t=0$ for each $t\in(0,r)$. Now we prove that
${\mathbf u}_j^t\in H_{0,\sigma}^1(D\setminus\overline\om)$ for
each $t\in (0,r)$. Indeed, for each fixed $t\in (0,r)$,
$spt{\mathbf u}_j^t=spt{\mathbf u}_j+t{\mathbf v}_j$ and $
{\mathbf u}_j^t=0$ over $\om+t{\mathbf v}_j$. Thus $spt{\mathbf
u}_j^t\subset\subset V_j\setminus\overline\om$. Therefore, there
exists an open Lipschitz subset $U_{j,t}\subset {\bf R}^n$ such
that $spt{\mathbf u}_j^t\subset\subset U_{j,t}\subset\subset
V_j\setminus\overline\om$. It is clear that ${\mathbf u}_j^t\in
H_0^1(U_{j,t})^n$ since $U_{j,t}$ is Lipschitz. This implies that
${\mathbf u}_j^t\in D_0^1(U_{j,t})$, for div${\mathbf u}_j^t=0$.
Thus ${\mathbf u}_j^t\in H_{0,\sigma}^1(U_{j,t})$( see Lemma
1.2.1, Chapter 3, $\cite{kn:[7]}$). Hence, ${\mathbf u}_j^t\in
H_{0,\sigma}^1(D\setminus\overline\om)$. On the other hand, one
can easily verify that
$$\|{\mathbf u}_j^t-{\mathbf
u}_j\|_{H^1(D\setminus\overline\om)}\rightarrow 0\ \ \ \mbox{as}\
\ t\rightarrow 0^+.$$ Thus we have proved that ${\mathbf u}_j\in
H_{0,\sigma}^1(D\setminus\overline\om)$.

 Since ${\mathbf u}=\sum_{j=0}^m{\mathbf u}_j$, it follows
that ${\mathbf u}\in H_{0,\sigma}^1(D\setminus \overline\om)$.
This thus  completes the proof of Theorem 4.5. $\endpf$

\vskip 3mm

{\bf Theorem 4.6 } Let $\om\subset {\bf R}^n$ for $n=2,3$  be of
class $\cal C$. Then
$$D_0^1(\om)=H_{0,\sigma}^1(\om).$$

\vskip 3mm

{\bf Proof of Theorem 4.6}\quad It suffices to show that
$D_0^1(\om)\subset H_{0,\sigma}^1(\om)$. Let ${\mathbf u}\in
D_0^1(\om)$ and let $D,D_1$ be bounded open Lipschitz subsets in
${\bf R}^n$ such that $\om\subset\subset D_1\subset\subset D$.
Extend the function ${\mathbf u}$ to   ${\bf R}^n$ by setting it
to be zero outside $\om$. From the argument in the proof of
Theorem 4.5, it is clear that  ${\mathbf u}\in D_0^1(D)$. Let
$\tilde\om\equiv D_1\setminus\overline\om$. Then
$\tilde\om\subset\subset D$ is of class $\cal C$ since every
bounded open Lipschitz set is of class $\cal C$, and ${\bf u}=0$
over $\tilde \om$.  Let $\tilde\om'$ be given as in Lemma 4.2 (
i.e. $\tilde\om'$ can be retracted to $\tilde\om$) such that
$\tilde\om'\subset\subset D$,  and Let $\{V_j\}_{j=1}^m$ be the
open covering of $\p\om$, given as in (2.5) and (2.6), such that
$\p\om\subset\bigcup_{j=1}^mV_j\subset\subset\tilde\om'$. Then by
Theorem 4.1' and Theorem 4.4', there are functions $\{{\mathbf
u}_j\}_{j=1}^m$ and a neighborhood $\hat U$ of $\p\om$ in ${\bf
R}^n$ with $\hat U\subset\subset \bigcup_{j=1}^mV_j$ such that
$${\mathbf u}=\sum_{j=1}^m {\mathbf u}_j\ \ \mbox{ over } \hat U$$
and for all $j=1,\cdots, m,\ {\mathbf u}_j\in D_0^1(V_j),\
{\mathbf u}_j=0$ over $\tilde\om$. From which we imply ${\mathbf
u}_j=0$ a.e. in $D\setminus\overline\om$ for all $j=1,\cdots,m$.

Let ${\mathbf u}_0={\mathbf u}-\sum_{j=1}^m{\mathbf u}_j$. Then
${\mathbf u}_0\in D_0^1(D)$ and ${\mathbf u}_0=0$ a.e. in
$(D\setminus\overline\om )\bigcup\hat U$, from which we get
${\mathbf u}_0=0$ a.e. in $D\setminus(\om\setminus\overline{\hat
U})$. Since $\om\setminus\overline{\hat U}\subset \subset \om$,
one can find a Lipschitz open subset $\hat \om$ such that
$\om\setminus\overline{\hat
U}\subset\subset\hat\om\subset\subset\om$. Because ${\mathbf
u}_0=0$ a.e. in $D\setminus\overline{\hat \om}$, we get
 by the trace theorem that ${\mathbf u}_0\in D_0^1(\hat
\om)$. Then by Lemma 1.2.1 in Chapter 3 in $\cite{kn:[7]}$, it
follows that ${\mathbf u}_0\in H_{0,\sigma}^1(\hat \om)\subset
H_{0,\sigma}^1(\om)$.

Next for each $j=1,\cdots,m$, by setting ${\mathbf
u}_j^t(x)={\mathbf u}_j(x+t{\mathbf v}_j)$ for $t>0$ small enough
and using the same argument as that in the proof of Theorem 4.5,
we  obtain that ${\mathbf u}_j\in H_{0,\sigma}^1(\om)$.

Hence, we proved that ${\mathbf u}=\sum_{j=0}^m{\mathbf u}_j\in
H_{0,\sigma}^1(\om)$.
This completes the proof of Theorem 4.6. $\endpf$

\vskip 3mm

{\bf Remark 4.7 }\quad After some slight modification in the proof
of Theorem 3.1, we can verify a similar decomposition result for
vector-valued function ${\mathbf u}\in L^2(D)^n$ with
$\mbox{div}{\mathbf u}=0$ in $D$. (Here $\mbox{div}{\mathbf u}$ is
in the sense of distribution.) To be more precise, we let
$\{U_j\}_{j=1}^m\subset D$ be  open subsets such that
$\bigcup_{j=1}^m U_j\subset\subset D$. Then for any open subset
$\hat U$ with $\hat U\subset\subset \bigcup_{j=1}^m U_j$ there are
functions $\{{\mathbf u}_j\}_{j=1}^m\subset L^2(D)^n$ such that
${\mathbf u}=\sum_{j=1}^m{\mathbf u}_j$ over $\hat U$, and for all
$j=1,\cdots,m, \ \mbox{div}{\mathbf u}_j=0$ in $D$ in the sense of
distribution, $spt{\mathbf u}_j\subset U_j$ and $\|{\mathbf
u}_j\|_{L^2(D)}\leq \mbox{C}\|{\mathbf u}\|_{L^2(D)}$.

Furthermore, the same localization property as stated in Theorem
4.1 and Theorem 4.2 for such a partition still holds. Thus we can
also obtain the following: If $\om$ is of class $\cal C$ in ${\bf
R}^n,\ n=2,3$, then
$$L_\sigma^2(\om)=\{{\mathbf u}\in L^2({\bf R}^n)^n;\ \mbox{div}{\mathbf u}=0\ \
\mbox{in}\ \ {\bf R}^n,\ {\mathbf u}=0 \ \ \mbox{a.e. in}\ \ {\bf
R}^n\setminus\om\},$$
 where
 $L_\sigma^2(\om)=\overline{ C_{0,\sigma}^\infty(\om)}^{\|\cdot\|_{L^2}}$,
 the completion of $C_{0,\sigma}^\infty(\om)$ in the norm
 of $L^2(\om)^n$.

\vskip 3mm

\def\theequation{5.\arabic{equation}}
\setcounter{equation}{0} \noindent {\large\bf5.\quad The existence
of a shape optimization problem for the stationary Navier-Stokes
equation }

\vskip 3mm

Let $D$ and $B$ be   open bounded Lipschitz subsets in ${\bf
R}^n,n=2,3$, such that $B\subset\subset D$.  Let $\cal O$ be the
family of certain open subsets $\om$ contained in $B$, which is
defined as in section 2. Consider for each $\om\in\cal O$ the
stationary Navier-Stokes equation on $D\setminus\overline\om$ :
\begin{eqnarray}
\begin{cases}{-\gamma\Delta {\mathbf u}-({\mathbf
u}\cdot\nabla){\mathbf u}+\nabla p={\mathbf f} &{\rm in} \
$D\setminus\overline\om$,\cr \mbox{div}{\mathbf u}=0 &{\rm in} \
$D\setminus\overline\om$,\cr u=0 &{\rm on} \ $\p D\bigcup\p\om$,}
\end{cases}
\end{eqnarray}
where ${\mathbf f}$ is a given function in $L^2(D)^n$. We shall
denote by ${\mathbf u}_\om$ the weak solution of (5.1), i.e.,
${\mathbf u}_\om\in H_{0,\sigma}^1(D\setminus\overline\om)$ and
$$\gamma\int_{D\setminus\overline\om}\nabla {\mathbf
u}_\om\cdot\nabla\Phi dx+\int_{D\setminus\overline\om}({\mathbf
u}_\om\cdot\nabla){\mathbf u}_\om\cdot\Phi
dx=\int_{D\setminus\overline\om}{\mathbf f}\cdot\Phi dx$$ for all
$\Phi\in C_{0,\sigma}^\infty(D\setminus\overline
\om)\equiv\{\Psi\in C_0^\infty(D\setminus\overline\om)^n;
\mbox{div}\Psi=0\}$.

It is well know that ( see $\cite{kn:[7]}$ and $\cite{kn:[4]}$ )
for each $\om\in\cal O$,  (5.1) has at least one weak solution.
Moreover, there is a positive constant C$(\gamma,d)$ depending
only on $\gamma$ and the width of $D$ such that if
\begin{eqnarray}
\|{\mathbf f}\|_{L^2(D)^n}<\mbox{C}(\gamma,d),
\end{eqnarray}
then the weak solution to (5.1) is unique. In the following we
always assume that (5.2) holds.

Our purpose in this section is  to study the following shape optimization problem:

 $ (P)\quad
\quad\;\;\;\;\;\;\;\;\;   \;\;\;\;\;\;\;\;\;\; \mbox{Min}_{\om\in
\cal O}\;\d{\int_{B\setminus \o \om}} J(x, {\bf u}_{\om}, \nabla
{\bf u}_{\om})\,dx\;\;\;\;\;\;\;\;\;\;\;\;\;\;\;\;\;\;\;\;\; $\\ subject
to Equation (5.1).

Here we assume that $J$ satisfies the following:

$(H)\quad\quad  J: D\times {\bf R}^n\times {\bf R}^{n\times
n}\rightarrow {\bf R}^+\equiv [0,\infty)$ is measurable in the
first variable and continuous in others. Moreover, there exist a
positive constant C and a function $g\in L^1(D)$ such that
$\mbox{for all} \quad (x,\xi,\eta)\in D\times {\bf R}^n\times {\bf
R}^{n\times n}$,
\begin{eqnarray}
J(x,\xi,\eta)\leq \mbox C (g(x)+|\xi|^2_{{\bf R}^n}+|\eta|^2_{{\bf R}^{n\times
n}}).\quad
\end{eqnarray}

Problem $(P)$ is to ask the shape of a body $\om$ in $\cal O$ such
that the cost functional takes its minimum. For instance, if we
take $J(x,{\mathbf u}_\om,\nabla {\mathbf u}_\om)\equiv |\nabla
{\mathbf u}_\om|^2+|(\nabla {\mathbf u}_\om)^T|^2$, where $(\nabla
{\mathbf u})^T$ denotes the transpose of $\nabla {\mathbf u}_\om$,
then problem $(P)$ is to ask the shape of a body $\om$ among $\cal
O$ having the smallest drag. ( see $\cite{kn:[6]}$ .)

\vskip 3mm

{\bf Theorem 5.1 } Assume that the conditions in (4.2) and (4.3)
are fulfilled. Then the shape optimization problem $(P)$ has at
least one solution $\om^*\in\cal O$.

\vskip 3mm

{\bf Proof of Theorem 5.1} Let $\{\om_m\}_{m=1}^\infty\subset \cal
O$ be a minimizing sequence for Problem $(P)$. By Lemma 2.1, there
exist some $\om^*\in \cal O$ and a subsequence of
$\{\om_m\}_{m=1}^\infty$, still denoted by itself, such that
$\mbox{Hlim}_{m\rightarrow\infty}\om_m=\om^*$. We notice that
$\om^*\subset\subset D$, for $B\subset\subset D$. Let ${\mathbf
u}_m,m=1,2,\cdots$, be the weak solutions to equation (4.1) with
$\om=\om_m$. Then ${\mathbf u}_m \in
H_{0,\sigma}^1(D\setminus\overline\om_m)$ and
$$\gamma \int_{D\setminus\overline\om_m}\nabla {\mathbf
u}_m\cdot\nabla\Phi dx+\int_{D\setminus\overline\om_m}({\mathbf
u}_m\cdot\nabla){\mathbf u}_m\cdot\Phi
dx=\int_{D\setminus\overline\om_m}{\mathbf f}\cdot\Phi dx$$ for
all $\Phi\in H_{0,\sigma}^1(D\setminus\overline\om_m)$. This
immediately gives that
$$\int_{D\setminus\overline\om_m}|\nabla {\mathbf
u}_m|^2 dx=    \int_{D\setminus\overline\om_m}{\mathbf
f}\cdot{\mathbf u}_m dx,$$ for
$\displaystyle\int_{D\setminus\overline\om_m}({\mathbf
u}_m\cdot\nabla){\mathbf u}_m\cdot{\mathbf u}_mdx=0$. Thus
$$\int_{D\setminus\overline\om_m}|\nabla {\mathbf
u}_m|^2 dx\leq \mbox C,$$ where C is a positive constant
independent of $m$.

Let
\begin{eqnarray*}
\hat {\mathbf u}_m(x)=\begin{cases}{{\mathbf u}_m(x) \ &{\rm in} \
$D\setminus\overline\om_m$,\cr 0\ &{\rm in}\ $\om_m\bigcup ({\bf
R}^n\setminus D)$.} \end{cases}
\end{eqnarray*}
Then it is clear that $\hat {\mathbf u}_m\in
H_{0,\sigma}^1(D\setminus\overline\om_m)\subset H_{0,\sigma}^1(D)$
and $\{\hat {\mathbf u}_m\}_{m=1}^\infty$ is bounded in
$H_0^1(D)^n$. Thus there exists a subsequence of $\{\hat {\mathbf
u}_m\}_{m=1}^\infty$, still denoted by itself, such that
\begin{eqnarray}
\hat {\mathbf u}_m\rightarrow\hat{\mathbf u}\ \ \mbox{weakly in }\
\ H^1(D)^n \ \ \mbox{and strongly in }\ \ L^2(D)^n
\end{eqnarray}
for some $\hat {\mathbf u}\in H_{0,\sigma}^1(D)$.

We claim that $\hat {\mathbf u}=0\ \ \mbox{a.e. in }\  \om^*$.
Indeed, for any open subset $K$ with $K\subset\subset\om^*$, it
follows from Lemma 2.2 that these exists an integer $m_K>0$ such
that $\overline K\subset\om_m$ for all $m\leq m_K$. Thus
$$\int_K|\hat {\mathbf
u}(x)|^2dx=\mbox{lim}_{m\rightarrow\infty}\int_K|\hat {\mathbf
u}_m(x)|^2dx\leq\overline{\mbox{lim}}_{m\rightarrow\infty}\int_{\om_m}|\hat
{\mathbf u}_m(x)|^2dx=0,$$
which implies $\hat {\mathbf u}(x)=0$
a.e. in $K$ and so $\hat {\mathbf u}=0$ a.e. in $\om^*$.

Now it follows from Theorem 4.5 that $\hat {\mathbf u}\in
H_{0,\sigma}^1(D\setminus\overline\om^*)$.

We next prove
\begin{eqnarray}
\gamma\int_{D\setminus\overline\om^*}\nabla\hat {\mathbf
u}\cdot\nabla\Phi dx+\int_{D\setminus\overline\om^*}(\hat {\mathbf
u}\cdot\nabla)\hat {\mathbf u}\cdot\Phi
dx=\int_{D\setminus\overline\om^*}{\mathbf f}\cdot\Phi dx
\end{eqnarray}
for any $\Phi\in C_{0,\sigma}^\infty(D\setminus\overline\om^*)$.
To this aim, we fix $\Phi\in
C_{0,\sigma}^\infty(D\setminus\overline\om^*)$ and let
$K=(spt\Phi)^o$, the interior of $spt\Phi$. Then by Lemma 2.3,
there is a positive integer $m(K)$ such that $\overline K\subset
D\setminus \om_m$ for any $m\geq m(K)$. Therefore,
$$\gamma\int_{spt\Phi}\nabla\hat {\mathbf
u}_m\cdot\nabla\Phi dx+\int_{spt\Phi}\hat{\mathbf u}_m\cdot\Phi
dx=\int_{spt\Phi}{\mathbf f}\cdot\Phi dx$$ for any $m\geq m(K)$.
Passing to the limit, as $m\rightarrow\infty $, in the above, we
get (5.5). Hence,  $\hat {\mathbf u}$ is a weak solution to
equation (5.1) where $\om=\om^*$.

Finally, we shall show
$$\int_{D\setminus\overline\om^*}J(x,\hat{\mathbf
u},\nabla \hat {\mathbf u})dx=\mbox{Min}_{\om\in\cal
O}\int_{D\setminus\overline\om}J(x,{\mathbf u}_\om,\nabla {\mathbf
u}_\om )dx.$$

For this purpose, we let
$D\setminus\overline\om^*=\bigcup_{j=1}^\infty G_j$ where
$G_j,j=1,\cdots$, are open subsets in $D\setminus\overline\om^*$
such that $\overline G_j\subset G_{j+1}$. Then by Lemma 2.3, we
obtain  that for each $j$, there exists an integer $m_j>0$ such
that $\overline G_j\subset D\setminus \overline\om_m $ for any
$m\geq m_j$. Therefore, we have  for each $j$, that
\begin{eqnarray}
\underline\lim_{m\rightarrow\infty}\int_{D\setminus\overline\om_m}J(x,\hat{\mathbf
u}_m,\nabla\hat {\mathbf u}_m )dx\geq
\underline\lim_{m\rightarrow\infty}\int_{G_j}J(x,\hat{\mathbf
u}_m,\nabla\hat{\mathbf u}_m) dx.
\end{eqnarray}
On the other hand, since $\hat {\mathbf u}$ is the weak solution
of (5.1) where $\om=\om^*$, we have
$$\gamma\int_{D\setminus\overline\om^*}|\nabla\hat{\mathbf
u}|^2dx=\int_{D\setminus\overline\om^*}{\mathbf
f}\cdot\hat{\mathbf u}dx,$$
 from which it follows that
\begin{eqnarray*}
\int_D|\nabla\hat{\mathbf u}_m|^2dx-\int_D|\nabla\hat{\mathbf
u}|^2dx &&=\int_{D\setminus\overline\om_m}|\nabla {\mathbf
u}_m|^2dx-\int_{D\setminus\overline\om^*}|\nabla\hat{\mathbf
u}|^2dx\\ &&=\f{1}{\gamma}(\int_{D\setminus\overline\om_m}{\mathbf
f}\cdot {\mathbf
u}_mdx-\int_{D\setminus\overline\om^*}{\mathbf f}\cdot\hat{\mathbf u}dx)\\
&&=\frac{1}{\gamma}\int_D{\mathbf f}\cdot\hat{\mathbf
u}_mdx-\int_D{\mathbf f}\cdot\hat{\mathbf u}dx\rightarrow 0 \ \
\mbox{as}\ \ m\rightarrow\infty,
\end{eqnarray*}
i.e.,
$$\|\hat{\mathbf u}_m\|_{H^1(D)^n}\rightarrow \|\hat{\mathbf
u}\|_{H^1(D)^n}.$$ This together with (5.4) implies that
$$\hat{\mathbf u}_m\rightarrow\hat{\mathbf u}\ \ \mbox{strongly in
 }\ H^1(D)^n.$$
Now by the assumption in $(H)$, Lemma 2.3 and by Fatou's lemma, we get
 for each $j$, the following:
$$\underline\lim_{m\rightarrow\infty}\int_{G_j}J(x,\hat{\mathbf
u}_m,\nabla\hat{\mathbf u}_m)dx\geq\int_{G_j}J(x,\hat{\mathbf
u},\nabla\hat{\mathbf u})dx.$$ This together with (5.6) shows that
 for each $j$ fixed,
\begin{eqnarray*}
\mbox{Min}_{\om\in\cal O}\int_{D\setminus\overline\om}J(x,{\mathbf
u}_\om,\nabla{\mathbf u}_\om)dx&& \geq
\underline\lim_{m\rightarrow\infty}\int_{G_j}J(x,\hat{\mathbf
u}_m,\nabla\hat{\mathbf u}_m)dx\\
&&\geq\int_{G_j}J(x,\hat{\mathbf u},\nabla\hat{\mathbf u})dx.
\end{eqnarray*}
Using Fatou's lemma again, we get
\begin{eqnarray*}\displaystyle
\mbox{Min}_{\om\in\cal O}\int_{D\setminus\overline\om}J(x,{\mathbf
u}_\om,\nabla{\mathbf u}_\om)dx&& \geq
\underline\lim_{j\rightarrow\infty}\int_{D\setminus\overline\om^*}\chi_{G_j}J(x,\hat{\mathbf
u} ,\nabla\hat{\mathbf u} )dx\\
&& \geq
\int_{D\setminus\overline\om^*}\underline\lim_{j\rightarrow\infty}\chi_{G_j}J(x,\hat{\mathbf
u} ,\nabla\hat{\mathbf u} )dx\\
&&=\int_{D\setminus\overline\om^*}J(x,\hat{\mathbf
u},\nabla\hat{\mathbf u})dx.
\end{eqnarray*}
where $\chi_{G_j}$ denotes the characteristic function of $G_j$.
Hence $\hat {\mathbf u}$ is an optimization solution for problem
$(P)$.

This completes the proof of Theorem 5.1. $\endpf$

\vskip 3mm

{\bf Remark 5.2 } By the same argument as  above and making use of
Theorem 4.6 instead, we can similarly obtain the existence for the
following shape optimization problem:

$\displaystyle\quad (\tilde P)\quad\quad\quad\quad\quad\quad
\mbox{Min}_{\om\in\cal O}\int_\om J(x,{\mathbf u}_\om,\nabla
{\mathbf u}_\om) dx$\\
subject to
\begin{eqnarray}
\cases{-\gamma\Delta {\mathbf u}-({\mathbf u}\cdot\nabla){\mathbf
u}+\nabla p={\mathbf f}\ &{\rm in}\ $\om$,\cr \mbox{div}{\mathbf
u}=0\ &{\rm in} \ $\om$,\cr {\mathbf u}=0\ &{\rm on}\ $\p\om$,}
\end{eqnarray}
where ${\mathbf u}_\om$ denotes the weak solution of (5.7)
corresponding to $\om\in \cal O$.

\end{document}